\documentclass[11pt]{amsart}
\usepackage{amsfonts,amsmath,amssymb,amsthm}
\usepackage{hyperref}
\usepackage[none]{hyphenat}
\usepackage{tikz}
\usetikzlibrary{shapes.arrows,decorations.pathreplacing} 

\usepackage{array}
\newcolumntype{R}[1]{>{\raggedleft\let\newline\\\arraybackslash\hspace{0pt}}m{#1}}

\theoremstyle{plain}
\newtheorem{theorem}{Theorem}[section]

\newtheorem{lemma}[theorem]{Lemma}

\theoremstyle{definition}
\newtheorem{definition}[theorem]{Definition}
\newtheorem*{example}{Example}
\newtheorem*{remark}{Remark}

\newcommand{\size}[1]{\lvert{#1}\rvert}

\newcommand{\gray}[1]{\textcolor{lightgray}{#1}}

\title[A natural bijection for contiguous pattern avoidance in words]{A natural bijection for contiguous\\ pattern avoidance in words}

\author{Julia Carrigan}
\address{Mathematics Department\\ Occidental College\\ Los Angeles, CA\\ USA}

\author{Isaiah Hollars}
\address{Department of Mathematics\\University of South Carolina\\Columbia, SC\\USA}

\author{Eric Rowland}
\address{
	Department of Mathematics \\
	Hofstra University \\
	Hempstead, NY \\
	USA
}

\thanks{This work was done in the 2022 New York City Discrete Mathematics REU and was supported by NSF award DMS-2051026.}

\begin{document}
\maketitle

\begin{abstract}
		Two words $p$ and $q$ are avoided by the same number of length-$n$ words, for all $n$, precisely when $p$ and $q$ have the same set of border lengths.
		Previous proofs of this theorem use generating functions but do not provide an explicit bijection. We give a bijective proof for all pairs $p, q$ that have the same set of proper borders, establishing a natural bijection from the set of words avoiding $p$ to the set of words avoiding $q$.
\end{abstract}

\section{Introduction}\label{intro}

    Combinatorialists have studied pattern avoidance in multiple contexts. In this paper, we are interested in the avoidance of contiguous patterns in words.
We say that a word $w$ \emph{avoids} a word $p$ if $w$ does not contain a contiguous occurrence of $p$. We refer to the word $p$ as a \emph{pattern}.
For example, the word $010$ avoids the pattern $00$ but does not avoid $10$.
Let $\mathbb{N}$ denote the set of non-negative integers.
	\begin{definition}
		Let $p$ and $q$ be two words on a finite alphabet $\Sigma$. Define \[A_n(p)=\{w \in \Sigma^n : \text{w avoids $p$}\}.\]
		The words $p$ and $q$ are \emph{avoidant-equivalent} if $\size{A_n(p)} = \size{A_n(q)}$ for all $n \in \mathbb{N}$.
	\end{definition}

    This notion of equivalence is analogous to Wilf equivalence for non-contiguous permutation patterns, which has been studied extensively.
	When two permutation patterns are avoided by the same number of permutations, researchers seek a bijective explanation.
	See for example the survey by Claesson and Kitaev~\cite{Claesson--Kitaev} of bijections between permutations that avoid $321$ and permutations that avoid $132$.

	Analogously, when two words $p, q$ are avoidant-equivalent, we would like a natural bijection from $A_n(p)$ to $A_n(q)$ for each $n$, since this provides a combinatorial explanation for the equivalence and therefore a deeper understanding of the relationship between these two structures.
    One can obtain a trivial bijection from $A_n(p)$ to $A_n(q)$ by first sorting the two sets lexicographically and then mapping the $i$th word of $A_n(p)$ to the $i$th word of $A_n(q)$. However, this bijection is mostly arbitrary; it assumes we already know that $\size{A_n(p)} = \size{A_n(q)}$, and it is computationally intensive, since computing the image of a particular word requires first computing the complete sets $A_n(p)$ and $A_n(q)$. In this paper, we establish a more natural bijection for many pairs of patterns. In particular, our bijection provides a new proof of avoidant-equivalence for these pairs.
    
    A sufficient condition for two patterns to be avoidant-equivalent has essentially been known since the work of Solov’ev~\cite{solovev2}, who determined the expected time required for a pattern $p$ to appear in a word built randomly letter by letter. Solov’ev showed that the expected time depends only on the lengths of the borders of $p$.

	\begin{definition}
		Let $p$ be a word. A non-empty word $x$ is a \emph{border} of $p$ if $x$ is both a prefix and a suffix of $p$.
		Let \[b(p)=\{\size{x}: \text{$x$ is a border of $p$} \}.\] We call $b(p)$ the \emph{border length set} of $p$. A non-empty word $x$ is a \emph{proper border} of $p$ if $x$ is a border of $p$ and $x \neq p$.
	\end{definition}

	\begin{example}
		Let $\Sigma=\{0,1\}$. The borders of $p = 0110$ are $0$ and $0110$. These borders can be thought of as the ways $p$ can overlap itself:
		\begin{center}
			\begin{tikzpicture}
				\node[] at (0,2) {0110};
				\node[] at (0.6,2.5) {0110};
				
				\node[] at (2,2) {0110};
				\node at (2.4,2.5) {0110};
				
				\node at (4,2) {0110};
				\node at (4.2,2.5) {0110};
				
				\node at (6,2) {0110};
				\node at (6,2.5) {0110};
				
				\draw[green] (0.18,1.77) rectangle (0.43,2.72);
				\draw [red] (2, 1.77) rectangle (2.4, 2.72);
				\draw [red] (3.8, 1.77) rectangle (4.39, 2.72);
				\draw [green] (5.61, 1.77) rectangle (6.39, 2.72);
			\end{tikzpicture}
		\end{center}
		The border length set is $b(0110)= \{1, 4\}$.
  The only proper border of $p$ is $0$.
	\end{example}

	It follows from the paper of Solov’ev, and more explicitly from the work of Guibas and Odlyzko~\cite{Guibas--Odlyzko}, that if $b(p) = b(q)$ then $p$ and $q$ are avoidant-equivalent.
	Moreover, Guibas and Odlyzko give a method for computing the generating function of the number of words avoiding a pattern (or set of patterns).
	Let $k = \size{\Sigma}$ be the size of the alphabet, let $l = \size{p}$, and define the polynomial $B(x)=\sum_{i \in b(p)}x^{l-i}$.
	Then
	\begin{equation}\label{gen func eq}
		\sum_{n\geq 0 }\size{A_n(p)}x^{n} = \frac{B(x)}{(1-kx)B(x)+x^{l} }.
	\end{equation}
	This generating function was obtained by Kim, Putcha, and Roush~\cite{Kim--Putcha--Roush} and Zeilberger~\cite{Zeilberger}. It can also be obtained by the Goulden--Jackson cluster method~\cite{Goulden--Jackson}; see the treatment by Noonan and Zeilberger~\cite{noonan} for a friendly introduction.

	\begin{example}
		For the word $p=0110$, the border length set is $b(p) =\{1,4\}$. For $q=1011$, we have $b(q)=\{1,4\}$ as well.
		Therefore $b(p) = b(q)$, and the series expansion of $\frac{x^3+1}{1-2x+x^3-x^4}$ gives the sizes of both $A_n(p)$ and $A_n(q)$ for all $n \in \mathbb{N}$. In particular, $p$ and $q$ are avoidant-equivalent.
	\end{example}
	
	The main result of this paper (Theorem~\ref{main theorem}) is the following.
	Suppose $p$ and $q$ are words on a finite alphabet $\Sigma$.
	If the set of proper borders of $p$ is equal to the set of proper borders of $q$, then the map $\phi_L$, which is defined in Section~\ref{preliminaries} and iteratively replaces occurrences of $q$ with $p$, is a bijection from $A_n(p)$ to $A_n(q)$ for all $n$.
	Note that here the condition is that the sets of proper borders themselves are equal, as opposed to the sets of border lengths.

	For words on the binary alphabet $\Sigma = \{0, 1\}$, there are $103764$ pairs of length-$10$ avoidant-equivalent patterns, and our theorem provides a bijection for $71058$ of these pairs.	
	Additionally, there are two types of trivial bijections --- left--right reversal and permutations of $\Sigma$.
	Compositions of all these bijections provide bijections for $103460$ pairs, which is $99.7\%$ of avoidant-equivalent pairs of length-$10$ patterns.
	See Section~\ref{how_many_section} and Table~\ref{bijections_table} for more data. The smallest pair of avoidant-equivalent patterns on $\{0, 1\}$ for which we do not have a natural bijection is $0010010$ and $0110110$, which have a border length set of $\{1,4,7\}$. 

\begin{example}
		Let $p=1001$ and $q=1101$. Since $p$ and $q$ have the same set of proper borders, namely $\{1\}$, the replacement function $\phi_L$ forms a bijection from $A_n(p)$ to $A_n(q)$. We would also like a bijection from $A_n(0110)$ to $A_n(q)$, since $b(0110) = \{1\} = b(q)$. The patterns $0110$ and $q$ do not have the same set of proper borders, since $0$ is a border of $0110$ but is not a border of $q$. However, if we let $\sigma$ be the letter permutation function, which replaces 0's with 1's and 1's with 0's, then $\sigma$ forms a bijection from $A_n(0110)$ to $A_n(p)$. Therefore the composition $\phi_L \circ \sigma$ is a bijection from $A_n(0110)$ to $A_n(q)$.
	\end{example}
	
	We mention that the sufficient condition $b(p) = b(q)$ for the patterns $p, q$ to be avoidant-equivalent is also necessary.
	This follows from the rational generating function in Equation~\eqref{gen func eq}, which provides a linear recurrence satisfied by $\size{A_n(p)}$. Namely, let $k=|\Sigma|$ and $l = \size{p}$ again, and let $s(n)= \size{A_n(p)}$. Then
	\[
		s(n) = k \, s(n-1) - s(n-l) + \sum_{\substack{i \in b(p) \\ i \neq l}}\Bigg(k \, s(n+i-l-1) - s(n+i-l) \Bigg).
	\]
	Using this recurrence, one can show that if $b(p) \neq b(q)$ then the sequence $(\size{A_n(p)})_{n\geq 0}$ first differs from $(\size{A_n(q)})_{n\geq 0}$ at 
	\[
	n=
	\begin{cases}
	   \min(\size{p},\size{q}) &\text{if } \size{p}\neq \size{q}\\
	   2\size{p}- \max\!\big(b(p)\,\triangle\, b(q)\big) &\text{if } \size{p}=\size{q} 
	\end{cases}
	\]
	where $\triangle$ denotes symmetric difference. Therefore, the patterns $p$ and $q$ are avoidant-equivalent if and only if $b(p)=b(q)$.  

	In Section~\ref{preliminaries}, we define replacement functions $\phi_L$ and $\phi_R$. Section~\ref{section: main proof} is dedicated to proving the main theorem, namely that $\phi_L$ establishes a bijection from $A_n(p)$ to $A_n(q)$ under the condition that the proper borders of $p$ and $q$ are identical.
	
	\section{Replacement functions}\label{preliminaries}
	In this section, we define the function $\phi_L$ that, under certain conditions, gives a bijection $\phi_L \colon A_n(p) \rightarrow A_n(q)$ in Section~\ref{section: main proof}.
	The general idea is to systematically replace each occurrence of $q$ in a word with $p$. We accomplish this with an iterative replacement process. We will define $\phi_L$ to take a $p$-avoiding word and scan from left to right looking for occurrences of $q$. If it finds $q$, it replaces the first occurrence of $q$ with $p$ and then starts the left-to-right scan over. The replacement process ends when no more $q$'s remain.
	We will prove in Lemma~\ref{lemma1} below that this process terminates.
	
	In the following definitions, we assume that we have two patterns $p$ and $q$ such that $b(p) = b(q)$. In particular, $\size{p} = \size{q}$.
	Let $f^{k}(w)$ be the word obtained by iteratively applying $k$ iterations of the function $f$ to $w$.

	\begin{definition}
		For a given $p$-avoiding word $w$, the \emph{single scan function $L$} replaces the leftmost $q$ in $w$ with $p$. If no $q$ exists, $L$ acts as the identity function.
		Define $\phi_L(w) = L^{i}(w)$, where $i$ is the least non-negative integer such that $L^{i}(w)$ contains no $q$'s. 
	\end{definition}

    Even though we are scanning left to right, a replacement in one position can be followed by a replacement to its left, as the following example shows.
	
	\begin{example}
		Let $p=011$ and $q=001$. The iterative replacement process of $\phi_L$ on the word $0001001 \in A_7(p)$ is as follows:
  \begin{alignat*}{4}
		    	 \gray{0}001\gray{001}& \stackrel{L}{\mapsto}\, &\gray{0}011\gray{001}&  &&&&\\
		    	&= &001\gray{1001}& \stackrel{L}{\mapsto}\, &011\gray{1001}& &&\\
		    	&&&= &\gray{0111}001& \stackrel{L}{\mapsto}\, &\gray{0111}011&\\
                &&&&& =&0111011&.
\end{alignat*}
		Thus, $\phi_L(0001001$) = 0111011. We have $0111011 \in A_7(q)$ as desired. 
	\end{example}
	
	To prove that $\phi_L$ forms a bijection from $A_n(p)$ to $A_n(q)$, we will prove that there exists a natural inverse function $\phi_R$. To this end, we define the functions $R$ and $\phi_R$, which are built to undo their counterparts $L$ and $\phi_L$. 
	
	\begin{definition}
		For a given $q$-avoiding word $w$, the \emph{single scan function $R$} replaces the rightmost $p$ in $w$ with $q$. If no $p$ exists, $R$ acts as the identity function.
	Define $\phi_R(w) = R^{j}(w)$, where $j$ is the least non-negative integer such that $R^{j}(w)$ contains no $p$'s. 
	\end{definition}
	
	\begin{example}
	    	Using $p = 011$ and $q = 001$ as in the previous example, one checks that $\phi_R(0111011) = 0001001$, so $\phi_R(\phi_L(0001001)) = 0001001$.
	\end{example}

	\begin{lemma}
		\label{lemma1}
		Let $p$ and $q$ be equal-length patterns such that $p \neq q$, and let $n \in \mathbb{N}$. For every $w \in A_n(p)$, we have $\phi_L(w) \in A_n(q)$. 
	\end{lemma}
	
	\begin{proof}
		Since $p \neq q$, either $p<q$ or $p>q$ lexicographically. Assume $p<q$, since the other case is analogous. If $w$ contains $q$, then $L(w) < w$.
		Therefore, iteratively applying $L$ produces lexicographically smaller words until the image no longer contains $q$. Since there are only finitely many  length-$n$ words on $\Sigma$, this happens after finitely many steps, at which point we have a word in $A_n(q)$.
	\end{proof}

	For a word $w$, we define $\overline{w}$ to be the reverse of $w$. Let $\overline{L}$ be the function that replaces the leftmost occurrence of $\overline{p}$ in a word with $\overline{q}$. Similarly, let $\overline{R}$ be the function that replaces the rightmost $\overline{q}$ in a word with $\overline{p}$.
	\begin{definition}
	    We now define the functions $\overline{\phi_L }\colon A_n(\Bar{q})\rightarrow A_n(\Bar{p})$ and $\overline{\phi_R }\colon A_n(\Bar{p})\rightarrow A_n(\Bar{q})$ in a similar fashion to $\phi_L$ and $\phi_R$. Define $\overline{\phi_L } = \overline{L}^{i}(w)$, where $i$ is the least non-negative integer such that $\overline{L}^{i}(w)$ contains no $\overline{p}$'s. Define $\overline{\phi_R } = \overline{R}^{j}(w)$, where $j$ is the least non-negative integer such that $\overline{R}^{j}(w)$ contains no $\overline{q}$'s. 
	\end{definition}
	
	\begin{lemma}
		\label{lemma2}
		Let $w \in A_n(p)$ and $v \in A_n(q)$. We have
		\begin{align} 
			\phi_R(v) &= \overline{\overline{\phi_L}(\overline{v})} \label{reverse 1}\\
			\phi_L(w) &= \overline{\overline{\phi_R}(\overline{w})}. \label{reverse 2}
		\end{align}
	\end{lemma}
	Intuitively, Equation~\eqref{reverse 2} says the functions $\phi_L$ and $\overline{\phi_R}$ are conjugate under word reversal. 
	\begin{example}
		Let $p=011$ and $q=001$, and let $w=0001001$. We will show Equation~\eqref{reverse 2} holds. An example above shows the computation of $\phi_L(w)=0111011$. Next we evaluate $\overline{\overline{\phi_R}(\overline{w})}$. Firstly, we have $\overline{w}=1001000$. Secondly, we evaluate $\overline{\phi_R}(\overline{w})$. Recall that $\overline{\phi_R}$ will scan right to left replacing $\overline{q}$'s with $\overline{p}$'s. The iterative replacement gives
\begin{alignat*}{4}
		    	\gray{100}100\gray{0}& \stackrel{\overline{R}}{\mapsto}\, &\gray{100}110\gray{0}& &&&&\\
		    	&={}&\gray{1001}100& \stackrel{\overline{R}}{\mapsto}\, &\gray{1001}110& &&\\
		    	 &&&={}&100\gray{1110}&  \stackrel{\overline{R}}{\mapsto}\, &110\gray{1110}&\\
         &&&&& =&1101110&.
\end{alignat*}
		This shows $\overline{\phi_R}(\overline{w})=1101110$. Since $\overline{1101110} = 0111011$, we have that $\overline{\overline{\phi_R}(\overline{w})} = 0111011$ as expected. 
	\end{example}
	
	\begin{proof}[Proof of Lemma~\ref{lemma2}]
	We prove Equation~\eqref{reverse 1} by induction on the number of replacement steps, denoted $k$. Then Equation~\eqref{reverse 2} will follow by symmetry.
 
 Let $j$ be the number of steps required by $\phi_R$ applied to $v$. We set out to show
	\begin{equation}
		\label{induc}
		R^k(v) = \overline{\big(\overline{L}^k(\overline{v})\big)},
	\end{equation}
	for $0\leq k \leq j$.
		It helps to first establish that, for any $v$ that still has some $p$ to replace, we have
		\begin{equation}
			\label{basic-conjug}
			R(v) = \overline{\big(\overline{L}(\overline{v})\big)}.
		\end{equation}
		To see why this is true, observe that replacing the rightmost $p$ is equivalent to
		\begin{itemize}
			\item reversing the word,
			\item replacing the leftmost $\Bar{p}$, and then
			\item reversing again.
		\end{itemize}

		For the base case, the left-hand side of Equation~\eqref{induc} equals $v$ because, when $k=0$, there are no $p$'s to replace in $v$. Similarly, $\overline{\big( \overline{L}^k(\overline{v}) \big) } = \overline{\big(\overline{v} \big)} =v $, because there are no $\Bar{p}$'s to replace in $\overline{v}$. 
		
		Inductively, assume Equation~\eqref{induc} holds for some value of $k$ where $0\leq k < j$.
		We have
		\begin{align*}
			\overline{\big(\overline{L}^{k+1}(\overline{v})\big)}&=
			\overline{\overline{L}\big(\overline{L}^{k}(\overline{v})\big)}\\
			&=  \overline{\overline{L}\big(\overline{R^k(v)}\big)}\text{\quad by the inductive hypothesis}\\
			&= R \big(R^k(v) \big) \text{\quad using Equation~\eqref{basic-conjug}}\\
			&= R^{k+1}(v).
		\end{align*}
		This establishes Equation~\eqref{induc}, which gives us Equation~\eqref{reverse 1}.
	\end{proof}

	\section{The main theorem}\label{section: main proof}
	
	With all this background, we are ready for the main result of the paper. 
	
	\begin{theorem}\label{main theorem}
		Let $\Sigma$ be a finite alphabet, and let $p$ and $q$ be distinct, equal-length words on $\Sigma$.
		If the set of proper borders of $p$ is equal to the set of proper borders of $q$, then $\phi_L\colon A_n(p) \rightarrow A_n(q)$ forms a bijection for all $n \in \mathbb{N}$.
	\end{theorem}
	For example, the set of proper borders for each of the words $0100$ and $0110$ is $\{0\}$. On the other hand, $0110$ and $1011$ do not have the same set of proper borders, despite $b(0110) = \{1,4\} = b(1011)$.

\begin{remark}

Let $w\in A_n(p)$. Observe that if $w$ also avoids $q$ then $\phi_L$ acts as the identity map on $w$. Therefore, Theorem~\ref{main theorem} implies that words that avoid $p$ and contain $q$ are in bijection with words that avoid $q$ and contain $p$.

A natural question is whether the number of $q$'s in $w$ is equal to the number of $p$'s in $\phi_L(w)$. While this is usually the case, there do exist counterexamples. For example, let $p=001$, $q=110$, and $w=110\gray{1}110$. After 3 replacements, we see that $\phi_L(w)= \gray{00}001\gray{01}$.
\end{remark}

	By Lemma~\ref{lemma1}, we have that $\phi_L$ is a map from $A_n(p)$ to $A_n(q)$. To prove Theorem~\ref{main theorem}, it suffices to show that $\phi_L$ is a bijection. To do this, we will show that $\phi_R$ is its inverse function, namely that $\phi_R \big(\phi_L(w) \big) = w $ for $w \in A_n(p)$ and also that $\phi_L \big(\phi_R(w) \big) = w $ for $w \in A_n(q)$. More specifically, we show that each one-step replacement $L$ that takes place in $\phi_L(w)$ is undone by a one-step replacement $R$.
	
	\begin{proof}[Proof of Theorem~\ref{main theorem}]
		Let $n \in \mathbb{N}$ and $w\in A_n(p)$. Let $i$ be the number of steps required by $\phi_L$ applied to $w$.
		We will show by induction that $L^{k-1}(w) = R\big(L^{k}(w)\big)$ for all $k$ satisfying $1 \leq k \leq i$, so that $R$ is the left inverse of $L$. It will then follow that $\phi_R\big(\phi_L(w) \big)=w$.
		Let $v \in A_n(q)$; then
		\begin{align*}
			\phi_L \big(\phi_R(v)\big)&= 
			\phi_L\big(\overline{\overline{\phi_L}(\overline{v}) } \big) \text{\quad by  Equation~\eqref{reverse 1}}\\
			&= \overline{ \overline{\phi_R} \big(\overline{\phi_L}(\overline{v})\big)} \text{\quad by Equation~\eqref{reverse 2}, letting $w =\big(\overline{\overline{\phi_L}(\overline{v}) } \big)$}\\
			&= \overline{\big(\overline{v}\big)}\text{\quad because $\overline{\phi_R}$ is the left inverse of $\overline{\phi_L}$}\\
			&= v.
		\end{align*}
		Thus, we will have also shown that $\phi_L\big(\phi_R(v) \big)=v$, so that $\phi_R$ is both the left inverse and right inverse of $\phi_L$. It will follow that $\phi_L\colon A_n(p) \rightarrow A_n(q)$ is a bijection.

		It remains to prove that $L^{k-1}(w) = R\big(L^{k}(w)\big)$.
		For the base case $k=1$, the left-hand side of $L^{k-1}(w) = R\big(L^{k}(w)\big)$ is equivalent to applying zero $L$ operations on $w$, so it trivially equals $w$. The right-hand side of this equation is $R\big(L(w)\big)$. We denote the new $p$ inserted by $L$ as $\hat{p}$. We claim that $\hat{p}$ is the rightmost $p$ in $L(w)$; then the $R$ step function will find it first and will replace $\hat{p}$ back with a $q$.
  
        To prove the claim, assume that $\hat{p}$ is not rightmost in $L(w)$. Then there is a $p$ to the right of $\hat{p}$. If this $p$ does not overlap $\hat{p}$, then it would have also been present in $w$. But $w$ is $p$-avoiding; therefore $p$ must overlap $\hat{p}$. We denote the overlap in $L(w)$ as $x$: 
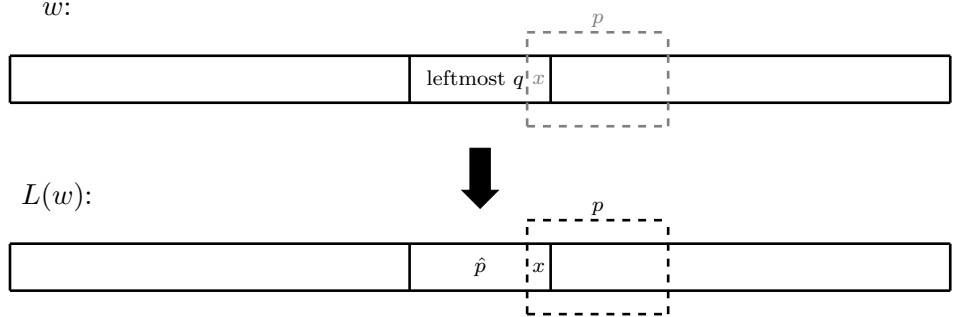
\begin{center}
\begin{tikzpicture}[scale=0.625]
	\node[] at (1,26) {$w$:};
			\draw[line width = 0.35mm] (0,24) -- (20,24);
			\draw[line width = 0.35mm] (0,25) -- (20,25);
			\draw[line width = 0.35mm] (0,24) -- (0,25);
			\draw[line width = 0.35mm] (20,25) -- (20,24);
			
			\draw[line width = 0.35mm] (8.5,24) -- (8.5,25);
			\draw[line width = 0.35mm] (11.5,24) -- (11.5,25);

			\draw[dashed,gray,opacity=0.5, line width = 0.35mm](11,23.5) -- (11,25.5);
			\draw[dashed,gray,opacity=0.5, line width = 0.35mm](14,23.5) -- (14,25.5);
			\draw[dashed,gray,opacity=0.5, line width = 0.35mm](11,23.5) -- (14,23.5);
			\draw[dashed,gray,opacity=0.5, line width = 0.35mm](11,25.5) -- (14,25.5);
			\node[font=\scriptsize] at (9.9,24.5) {leftmost $q$};
			\node[color=gray,opacity=0.5,font=\scriptsize] at (11.25,24.5) {$x$};
			\node[color=gray,opacity=0.5,font=\scriptsize] at (12.5,25.75) {$p$};
					
			\node[single arrow, draw=black, fill=black,
			minimum width = 10pt, single arrow head extend=3pt,
			minimum height=8mm,rotate=-90] at (10,22.5) {};

	\node[] at (1,22) {$L(w)$:};
			\draw[line width = 0.35mm] (0,20) -- (20,20);
			\draw[line width = 0.35mm] (0,21) -- (20,21);
			\draw[line width = 0.35mm] (0,20) -- (0,21);
			\draw[line width = 0.35mm] (20,21) -- (20,20);
			
			\draw[line width = 0.35mm] (8.5,20) -- (8.5,21);
			\node[font=\scriptsize] at (10,20.5) {$\hat{p}$};
			\draw[line width = 0.35mm] (11.5,20) -- (11.5,21);
			
			\draw[dashed, line width = 0.35mm](11,19.5) -- (11,21.5);
			\draw[dashed, line width = 0.35mm](14,19.5) -- (14,21.5);
			\draw[dashed, line width = 0.35mm](11,19.5) -- (14,19.5);
			\draw[dashed, line width = 0.35mm](11,21.5) -- (14,21.5);
			\node[font=\scriptsize] at (11.25,20.5) {$x$};
			\node[font=\scriptsize] at (12.5,21.75) {$p$};
\end{tikzpicture}
\end{center}
  The overlap $x$ is a border of $p$. Observe that since $p$ and $q$ have the same borders, the border segment $x$ is also in $w$ as a suffix of $q$. This means that $x$ wasn't altered when we swapped in $\hat{p}$. This implies that the overlapping $p$ is also in $w$. This contradicts our assumption that $w$ is $p$-avoiding. Therefore, $\hat{p}$ is the rightmost $p$ in $L(w)$, implying that $L^{k-1}(w) = R\big(L^{k}(w)\big)$ holds for $k = 1$.
		
		Inductively, assume that $L^{k-2}(w) = R\big(L^{k-1}(w)\big)$ for some $k$ between 1 and the number of steps required by $\phi_L$. This assumption means that once we replace the leftmost $q$ in $L^{k-2}(w)$ with $p$, this new $p$ must be the rightmost $p$ in $L^{k-1}(w)$ because we assumed that the $R$ function maps $L^{k-1}(w)$ back to $L^{k-2}(w)$ (the $R$ function scans from right to left); this is indicated by an arrow in each direction in the diagram below. To show the inductive hypothesis holds for $k+1$, we need to show this same relationship holds between words $L^{k-1}(w)$ and $L^{k}(w)$. Thus, we wish to show that, once we replace the leftmost $q$ in $L^{k-1}(w)$ with $p$, this new $p$ in $L^{k}(w)$ is the rightmost. The proof is split into four cases based on the possible positions of the leftmost $q$ in the word $L^{k-1}(w)$.
		
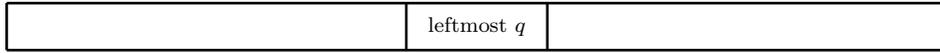
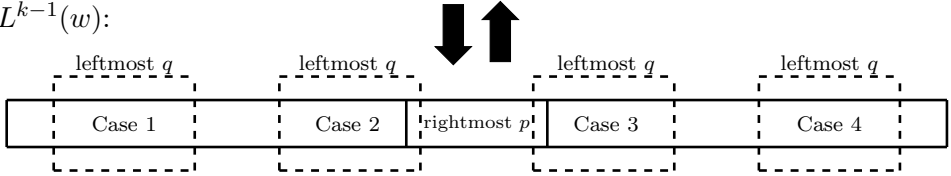
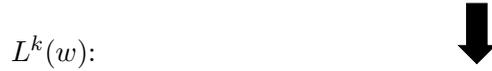
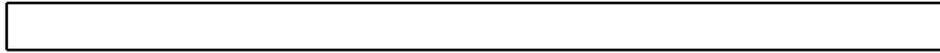
\begin{center}
\begin{tikzpicture}[scale=0.625]
\node[] at (1,30) {$L^{k-2}(w)$:};
			\draw[line width = 0.35mm] (0,28) -- (20,28);
			\draw[line width = 0.35mm] (0,29) -- (20,29);
			\draw[line width = 0.35mm] (0,28) -- (0,29);
			\draw[line width = 0.35mm] (20,29) -- (20,28);

			\draw[line width = 0.35mm] (8.5,28) -- (8.5,29);
			\node[font=\scriptsize] at (10,28.5) {leftmost $q$};
			\draw[line width = 0.35mm] (11.5,28) -- (11.5,29);

\node[single arrow, draw=black, fill=black,
minimum width = 10pt, single arrow head extend=3pt,
minimum height=8mm,rotate=-90] at (9.5,26.5) {}; 	
\node[single arrow, draw=black, fill=black,
minimum width = 10pt, single arrow head extend=3pt,
minimum height=8mm,rotate=90] at (10.5,26.35) {};

\node[] at (1,26.75) {$L^{k-1}(w)$:};
			\draw[line width = 0.35mm] (0,24) -- (20,24);
			\draw[line width = 0.35mm] (0,25) -- (20,25);
			\draw[line width = 0.35mm] (0,24) -- (0,25);
			\draw[line width = 0.35mm] (20,25) -- (20,24);
			
			\draw[line width = 0.35mm] (8.5,24) -- (8.5,25);
			\node[font=\tiny] at (10,24.5) {rightmost $p$};
			\draw[line width = 0.35mm] (11.5,24) -- (11.5,25);

			\draw[dashed, line width = 0.35mm](1,23.5) -- (1,25.5);
			\draw[dashed, line width = 0.35mm](4,23.5) -- (4,25.5);
			\draw[dashed, line width = 0.35mm](1,23.5) -- (4,23.5);
			\draw[dashed, line width = 0.35mm](1,25.5) -- (4,25.5);
			\node[font=\scriptsize] at (2.5,24.5) {Case 1};
			\node[font=\scriptsize] at (2.5,25.75) {leftmost $q$};
			
			\draw[dashed, line width = 0.35mm](5.8,23.5) -- (5.8,25.5);
			\draw[dashed, line width = 0.35mm](8.8,23.5) -- (8.8,25.5);
			\draw[dashed, line width = 0.35mm](5.8,23.5) -- (8.8,23.5);
			\draw[dashed, line width = 0.35mm](5.8,25.5) -- (8.8,25.5);
			\node[font=\scriptsize] at (7.25,24.5) {Case 2};
			\node[font=\scriptsize] at (7.25,25.75) {leftmost $q$};
			
			\draw[dashed, line width = 0.35mm](11.2,23.5) -- (11.2,25.5);
			\draw[dashed, line width = 0.35mm](14.2,23.5) -- (14.2,25.5);
			\draw[dashed, line width = 0.35mm](11.2,23.5) -- (14.2,23.5);
			\draw[dashed, line width = 0.35mm](11.2,25.5) -- (14.2,25.5);
			\node[font=\scriptsize] at (12.75,24.5) {Case 3 };
			\node[font=\scriptsize] at (12.75,25.75) {leftmost $q$};
			\draw[dashed, line width = 0.35mm](16,23.5) -- (16,25.5);
			\draw[dashed, line width = 0.35mm](19,23.5) -- (19,25.5);
			\draw[dashed, line width = 0.35mm](16,23.5) -- (19,23.5);
			\draw[dashed, line width = 0.35mm](16,25.5) -- (19,25.5);
			\node[font=\scriptsize] at (17.5,24.5) {Case 4};
			\node[font=\scriptsize] at (17.5,25.75) {leftmost $q$};

\node[single arrow, draw=black, fill=black,
minimum width = 10pt, single arrow head extend=3pt,
minimum height=8mm,rotate=-90] at (10,22.5) {};

\node[] at (1,22) {$L^{k}(w)$:};
			\draw[line width = 0.35mm] (0,20) -- (20,20);
			\draw[line width = 0.35mm] (0,21) -- (20,21);
			\draw[line width = 0.35mm] (0,20) -- (0,21);
			\draw[line width = 0.35mm] (20,21) -- (20,20);
			
\end{tikzpicture}
\end{center}

		\noindent
		\textbf{Case~1.}
		This position of $q$ in $L^{k-1}(w)$ implies that there is a $q$ in the same position in $L^{k-2}(w)$. But then we have a $q$ to the left of the leftmost $q$ in $L^{k-2}(w)$, a contradiction.
		
		\noindent
\textbf{Case~2. }Let $x$ be the overlap of the leftmost $q$ and the rightmost $p$ in $L^{k-1}(w)$. Suppose first that $x$ is a border of $q$. We use a similar argument as in the base case. Since $x$ did not change from $L^{k-2}(w)$ to $L^{k-1}(w)$, there must exist a $q$ in the same spot in $L^{k-2}(w)$. This $q$ is left of the leftmost $q$ in $L^{k-2}(w)$, so we have a contradiction.
  
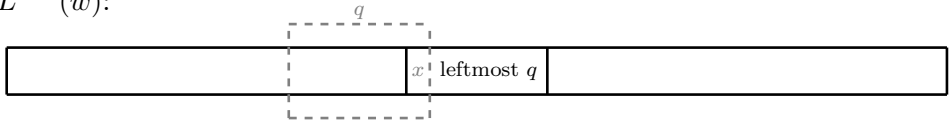
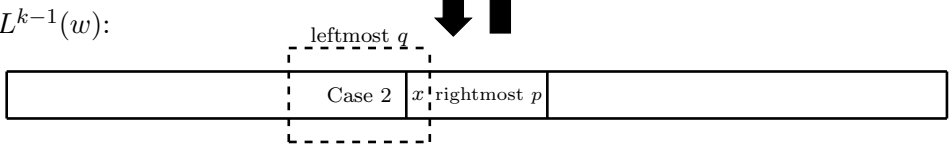
\begin{center}
\begin{tikzpicture}[scale=0.625]
	
\node[] at (1,26) {$L^{k-2}(w)$: };
	
			\draw[line width = 0.35mm] (0,24) -- (20,24);
			\draw[line width = 0.35mm] (0,25) -- (20,25);
			\draw[line width = 0.35mm] (0,24) -- (0,25);
			\draw[line width = 0.35mm] (20,25) -- (20,24);
			
			\draw[line width = 0.35mm] (8.5,24) -- (8.5,25);
			\node[font=\scriptsize] at (10.25,24.5) {leftmost $q$};
			\draw[line width = 0.35mm] (11.5,24) -- (11.5,25);
			
			\draw[dashed,gray,opacity=0.5, line width = 0.35mm](6,23.5) -- (6,25.5);
			\draw[dashed,gray,opacity=0.5, line width = 0.35mm](9,23.5) -- (9,25.5);
			\draw[dashed,gray,opacity=0.5, line width = 0.35mm](6,23.5) -- (9,23.5);
			\draw[dashed,gray,opacity=0.5, line width = 0.35mm](6,25.5) -- (9,25.5);
			\node[gray,opacity=0.5,font=\scriptsize] at (7.5,25.75) { $q$};
			\node[color=gray,opacity=0.5,font=\scriptsize] at (8.75,24.5) {$x$};

	\node[single arrow, draw=black, fill=black,
	minimum width = 10pt, single arrow head extend=3pt,
	minimum height=8mm,rotate=-90] at (9.5,22.5) {}; 	
	\node[single arrow, draw=black, fill=black,
	minimum width = 10pt, single arrow head extend=3pt,
	minimum height=8mm,rotate=90] at (10.5,22.35) {};

\node[] at (1,22) {$L^{k-1}(w)$: };
			\draw[line width = 0.35mm] (0,20) -- (20,20);
			\draw[line width = 0.35mm] (0,21) -- (20,21);
			\draw[line width = 0.35mm] (0,20) -- (0,21);
			\draw[line width = 0.35mm] (20,21) -- (20,20);
			
			\draw[line width = 0.35mm] (8.5,20) -- (8.5,21);
			\node[font=\tiny] at (10.25,20.5) {rightmost $p$};
			\draw[line width = 0.35mm] (11.5,20) -- (11.5,21);

			\draw[dashed, line width = 0.35mm](6,19.5) -- (6,21.5);
			\draw[dashed, line width = 0.35mm](9,19.5) -- (9,21.5);
			\draw[dashed, line width = 0.35mm](6,19.5) -- (9,19.5);
			\draw[dashed, line width = 0.35mm](6,21.5) -- (9,21.5);
			\node[font=\scriptsize] at (7.5,20.5) {Case 2};
			\node[font=\scriptsize] at (7.5,21.75) {leftmost $q$};
			\node[font=\scriptsize] at (8.75,20.5) {$x$};
	
\end{tikzpicture}
\end{center}

Suppose instead that $x$ is not a border of $q$ and therefore not a border of $p$. Toward a contradiction, assume that $\hat{p}$ is not the rightmost $p$ in $L^k(w)$. For this to occur, $\hat{p}$ must be overlapping with another $p$ on its right. We label this overlap segment $b$. Observe that $b$ is a border of $p$ and $q$. In particular, $b$ is a suffix of the leftmost $q$ and a prefix of the rightmost $p$ in $L^{k-1}(w)$. If $\size{b}> \size{x}$,  it would follow that $x$ is both a prefix and suffix of $b$. This would imply $x$ is a border of $p$, a contradiction.
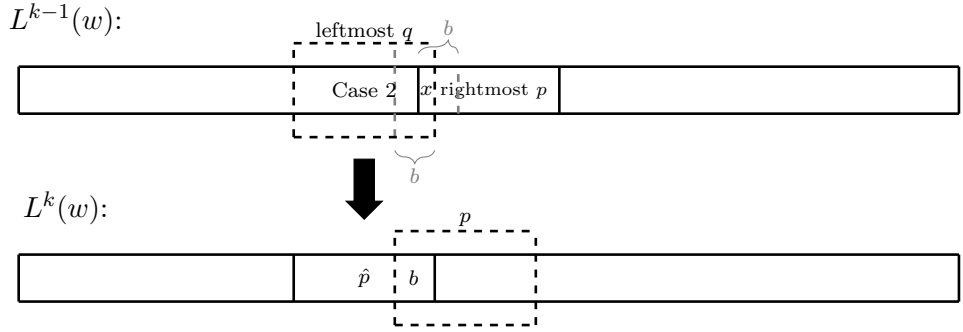
\begin{center}	
\begin{tikzpicture}[scale=0.625]
\node[] at (1,26) {$L^{k-1}(w)$: };	

			\draw[line width = 0.35mm] (0,24) -- (20,24);
			\draw[line width = 0.35mm] (0,25) -- (20,25);
			\draw[line width = 0.35mm] (0,24) -- (0,25);
			\draw[line width = 0.35mm] (20,25) -- (20,24);
			
			\draw[line width = 0.35mm] (8.5,24) -- (8.5,25);
			\node[font=\tiny] at (10.1,24.5) {rightmost $p$};
			\draw[line width = 0.35mm] (11.5,24) -- (11.5,25);
			\draw[dashed,line width = 0.35mm,color=gray,opacity=0.25] (9.35,24) -- (9.35,25);
			
			\draw [gray,opacity=0.5,decorate, decoration = {brace, mirror,raise=5pt, amplitude=5pt,aspect=0.5}] (8,23.5) -- (8.85,23.5)
			node[pos=0.5,below=8pt,font=\scriptsize]{$b$};
			
			\draw [gray,opacity=0.5,decorate, decoration = {brace,raise=5pt, amplitude=5pt,aspect=0.75}] (8.5,25) -- (9.35,25)
			node[pos=0.75,above=8pt,font=\scriptsize]{$b$};

			\draw[dashed, line width = 0.35mm](5.85,23.5) -- (5.85,25.5);
			\draw[dashed, line width = 0.35mm](8.85,23.5) -- (8.85,25.5);
			\draw[dashed, line width = 0.35mm,color=gray,opacity=0.25](8,23.5) -- (8,25.5);
			
			\draw[dashed, line width = 0.35mm](5.85,23.5) -- (8.85,23.5);
			\draw[dashed, line width = 0.35mm](5.85,25.5) -- (8.85,25.5);
			\node[font=\scriptsize] at (7.35,24.5) {Case 2};
			\node[font=\scriptsize] at (7.35,25.75) {leftmost $q$};
			\node[font=\scriptsize] at (8.675,24.5) {$x$};

			\node[single arrow, draw=black, fill=black,
			minimum width = 10pt, single arrow head extend=3pt,
			minimum height=8mm,rotate=-90] at (7.35,22.5) {}; 
	
\node[] at (1,22) {$L^{k}(w)$: };

			\draw[line width = 0.35mm] (0,20) -- (20,20);
			\draw[line width = 0.35mm] (0,21) -- (20,21);
			\draw[line width = 0.35mm] (0,20) -- (0,21);
			\draw[line width = 0.35mm] (20,21) -- (20,20);
			
			\draw[line width = 0.35mm] (5.85,20) -- (5.85,21);
			\node[font=\scriptsize] at (7.35,20.5) {$\hat{p}$};
			\draw[line width = 0.35mm] (8.85,20) -- (8.85,21);

			\draw[dashed, line width = 0.35mm](8,19.5) -- (8,21.5);
			\draw[dashed, line width = 0.35mm](11,19.5) -- (11,21.5);
			\draw[dashed, line width = 0.35mm](8,19.5) -- (11,19.5);
			\draw[dashed, line width = 0.35mm](8,21.5) -- (11,21.5);
			\node[font=\scriptsize] at (9.5,21.75) {$p$};
			\node[font=\scriptsize] at (8.4,20.5) {$b$};
			
\end{tikzpicture}
\end{center}
So it must be that $\size{b}< \size{x}$. Notice that the overlapping $p$ in $L^k(w)$ must have existed in the same position in $L^{k-1}(w)$ since $b$ was left unchanged when $\hat{p}$ was swapped in. This puts a $p$ to the right of the rightmost $p$ in $L^{k-1}(w)$, another contradiction.

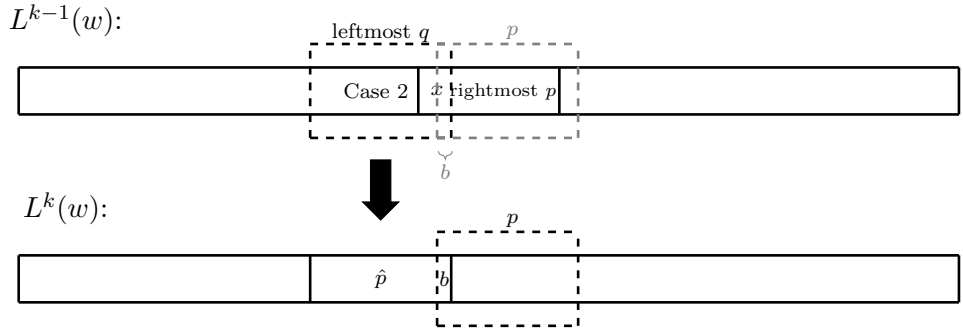
\begin{center}
\begin{tikzpicture}[scale=0.625]

\node[] at (1,26) {$L^{k-1}(w)$: };	

\draw[line width = 0.35mm] (0,24) -- (20,24);
\draw[line width = 0.35mm] (0,25) -- (20,25);
\draw[line width = 0.35mm] (0,24) -- (0,25);
\draw[line width = 0.35mm] (20,25) -- (20,24);
			
			\draw[line width = 0.35mm] (8.5,24) -- (8.5,25);
			\node[font=\tiny] at (10.3,24.5) {rightmost $p$};
			\draw[line width = 0.35mm] (11.5,24) -- (11.5,25);

			\draw[dashed, line width = 0.35mm](6.2,23.5) -- (6.2,25.5);
			\draw[dashed, line width = 0.35mm](9.2,23.5) -- (9.2,25.5);
			\draw[dashed, line width = 0.35mm](6.2,23.5) -- (9.2,23.5);
			\draw[dashed, line width = 0.35mm](6.2,25.5) -- (9.2,25.5);
			\node[font=\scriptsize] at (7.6,24.5) {Case 2};
			\node[font=\scriptsize] at (7.7,25.75) {leftmost $q$};
			\node[font=\scriptsize] at (8.9,24.5) {$x$};
			
			\draw [gray,opacity=0.5,decorate, decoration = {brace, mirror,raise=5pt, amplitude=3pt,aspect=0.5}] (8.925,23.5) -- (9.225,23.5)
			node[pos=0.5,below=6pt,font=\scriptsize]{$b$};
			
			
			\draw[dashed,gray,opacity=0.5, line width = 0.35mm](8.9,23.5) -- (8.9,25.5);
			\draw[dashed,gray,opacity=0.5, line width = 0.35mm](11.9,23.5) -- (11.9,25.5);
			\draw[dashed,gray,opacity=0.5, line width = 0.35mm](8.9,23.5) -- (11.9,23.5);
			\draw[dashed,gray,opacity=0.5, line width = 0.35mm](8.9,25.5) -- (11.9,25.5);
			\node[font=\scriptsize,gray,opacity=0.5] at (10.5,25.75) {$p$};
			
			\node[single arrow, draw=black, fill=black,
			minimum width = 10pt, single arrow head extend=3pt,
			minimum height=8mm,rotate=-90] at (7.7,22.5) {}; 
	
\node[] at (1,22) {$L^{k}(w)$: };	

			\draw[line width = 0.35mm] (0,20) -- (20,20);
			\draw[line width = 0.35mm] (0,21) -- (20,21);
			\draw[line width = 0.35mm] (0,20) -- (0,21);
			\draw[line width = 0.35mm] (20,21) -- (20,20);
			
			\draw[line width = 0.35mm] (6.2,20) -- (6.2,21);
			\node[font=\scriptsize] at (7.7,20.5) {$\hat{p}$};
			\draw[line width = 0.35mm] (9.2,20) -- (9.2,21);

			\draw[dashed, line width = 0.35mm](8.9,19.5) -- (8.9,21.5);
			\draw[dashed, line width = 0.35mm](11.9,19.5) -- (11.9,21.5);
			\draw[dashed, line width = 0.35mm](8.9,19.5) -- (11.9,19.5);
			\draw[dashed, line width = 0.35mm](8.9,21.5) -- (11.9,21.5);
			\node[font=\scriptsize] at (10.5,21.75) {$p$};
			\node[font=\scriptsize] at (9.05,20.5) {$b$};
			
\end{tikzpicture}
\end{center}

\noindent
\textbf{Case 3. } Suppose for contradiction that $\hat{p}$ is not the rightmost $p$ in $L^k(w)$. Then $\hat{p}$ must overlap with another $p$ to its right. We label the overlap $x$ in the diagram below. Since $x$ is a border of $p$, it is a border of $q$. Hence $x$ was left unchanged when $\hat{p}$ was substituted in. This implies that the $p$ to the right of $\hat{p}$ must have existed in the previous word $L^{k-1}(w)$. But this puts a $p$ to the right of the rightmost $p$ in $L^{k-1}(w)$, a contradiction.

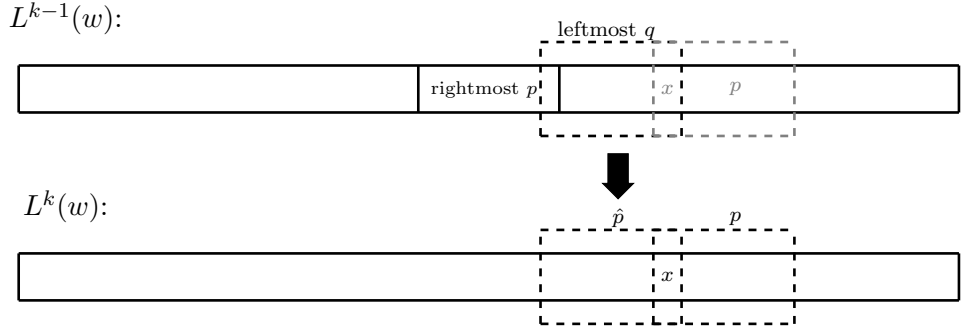
\begin{center}
\begin{tikzpicture}[scale=0.625]
\node[] at (1,26) {$L^{k-1}(w)$: };	
		
			\draw[line width = 0.35mm] (0,24) -- (20,24);
			\draw[line width = 0.35mm] (0,25) -- (20,25);
			\draw[line width = 0.35mm] (0,24) -- (0,25);
			\draw[line width = 0.35mm] (20,25) -- (20,24);
			
			\draw[line width = 0.35mm] (8.5,24) -- (8.5,25);
			\node[font=\tiny] at (9.9,24.5) {rightmost $p$};
			\draw[line width = 0.35mm] (11.5,24) -- (11.5,25);

			
			\draw[dashed, line width = 0.35mm](11.1,23.5) -- (11.1,25.5);
			\draw[dashed, line width = 0.35mm](14.1,23.5) -- (14.1,25.5);
			\draw[dashed, line width = 0.35mm](11.1,23.5) -- (14.1,23.5);
			\draw[dashed, line width = 0.35mm](11.1,25.5) -- (14.1,25.5);
			
			\node[font=\scriptsize] at (12.5,25.75) {leftmost $q$};
			\node[gray,opacity=0.5,font=\scriptsize] at (13.8,24.5) {$x$};
			
			\draw[dashed,gray,opacity=0.5, line width = 0.35mm](13.5,23.5) -- (13.5,25.5);
			\draw[dashed,gray,opacity=0.5, line width = 0.35mm](16.5,23.5) -- (16.5,25.5);
			\draw[dashed,gray,opacity=0.5, line width = 0.35mm](13.5,23.5) -- (16.5,23.5);
			\draw[dashed,gray,opacity=0.5, line width = 0.35mm](13.5,25.5) -- (16.5,25.5);
			\node[gray,opacity=0.5,font=\scriptsize] at (15.25,24.5){$p$};

			\node[single arrow, draw=black, fill=black,
			minimum width = 10pt, single arrow head extend=2pt,
			minimum height=6mm,rotate=-90] at (12.75,22.75) {};

\node[] at (1,22) {$L^{k}(w)$: };	
		
			\draw[line width = 0.35mm] (0,20) -- (20,20);
			\draw[line width = 0.35mm] (0,21) -- (20,21);
			\draw[line width = 0.35mm] (0,20) -- (0,21);
			\draw[line width = 0.35mm] (20,21) -- (20,20);
			
			\draw[dashed, line width = 0.35mm](11.1,19.5) -- (11.1,21.5);
			\draw[dashed, line width = 0.35mm](14.1,19.5) -- (14.1,21.5);
			\draw[dashed, line width = 0.35mm](11.1,19.5) -- (14.1,19.5);
			\draw[dashed, line width = 0.35mm](11.1,21.5) -- (14.1,21.5);
			\node[font=\scriptsize] at (12.75,21.75) { $\hat{p}$};

			\draw[dashed, line width = 0.35mm](13.5,19.5) -- (13.5,21.5);
			\draw[dashed, line width = 0.35mm](16.5,19.5) -- (16.5,21.5);
			\draw[dashed, line width = 0.35mm](13.5,19.5) -- (16.5,19.5);
			\draw[dashed, line width = 0.35mm](13.5,21.5) -- (16.5,21.5);
			\node[font=\scriptsize] at (15.25,21.75) {$p$};
			\node[font=\scriptsize] at (13.8,20.5) {$x$};

\end{tikzpicture}
\end{center}
  
	\noindent
	\textbf{Case 4.} Case 4 follows using the same argument as in Case~3. 
	\end{proof}
 
We contextualize the proof with an example and a counterexample.
\begin{example}
	Let $p = 0110$ and $q = 0010$. Note that the set of proper borders for both $p$ and $q$ is $\{0\}$, so $\phi_L$ is a bijection from $A_n(p)$ to $A_n(q)$. Let $w= 1001001011 \in A_{10}(p)$. This example demonstrates how each single scan function $L$ is undone by the function $R$. Observe that the first replacement aligns with the base case of the proof for Theorem \ref{main theorem}, while the second replacement aligns with Case 3B. Running $\phi_L$ on $w$ gives
        \begin{alignat*}{3}
				 \gray{1}0010\gray{01011}& \stackrel{L}{\mapsto}\, &\gray{1}0110\gray{01011}& && \\
				&={} &\gray{1011}0010\gray{11}& \stackrel{L}{\mapsto}\, &\gray{1011}0110\gray{11}&\\
                &&&= &1011011011&.
	\end{alignat*}
		Now we will run $\phi_L(w)=1011011011$ through $\phi_R$ to see that we get $w$ back. We also see that single scan $R$ successfully undoes every replacement made by an $L$. This gives us
  \begin{alignat*}{3}
				 \gray{1}0110\gray{11011}& \stackrel{R}{\mapsto}\, &\gray{1}0010\gray{11011}& &&\\
				&={}&\gray{1001}0110\gray{11}& \stackrel{R}{\mapsto}\, &\gray{1001}0010\gray{11}& \\
                &&&= &1001001011& = w .
\end{alignat*}
	\end{example} 
	
	\begin{example}
We now present a short counterexample. Let $p = 1011$ and $q = 0100$. Note that $b(p)=\{1,4\}=b(q)$, but $1$ is a proper border of $p$ and not a proper border of $q$. So, Theorem \ref{main theorem} does not guarantee $\phi_L$ will form a bijection from $A_n(p)$ to $A_n(q)$.  For the word $w_1 = 0101011 \in A_7(p)$, we have
        \begin{equation*}
				\gray{010}1011 \stackrel{L}{\mapsto} \gray{010}0100 = 0100100.
		\end{equation*}
For another word $w_2 = 1011100 \in A_7(p)$, we have
        \begin{equation*}
				1011\gray{100} \stackrel{L}{\mapsto} 0100\gray{100} = 0100100.
		\end{equation*}
Observe that $\phi_L(w_1) = 0100100 = \phi_L(w_2)$ so that $\phi_L$ does not provide a bijection.
\end{example}

	\subsection{How many bijections do we obtain?}\label{how_many_section}
    One might wonder how many of the possible bijections $\phi_L$ provides. We know that $\phi_L$ forms a bijection from $A_n(p)$ to $A_n(q)$ if $b(p) = \{\size{p}\} = b(q)$. Words with no proper borders, such as these, are known as \emph{borderless} words. The density of borderless words on a finite alphabet has been analyzed in detail.  Silberger~\cite{silberger} first discovered a recursive formula to count borderless words, and Holub \& Shallit~\cite{holub} investigated the probability that a random word is borderless. Notably, a long binary word $p$ chosen randomly has $\approx 27\%$ chance of being borderless and $\approx 30\%$ chance of having the border length set $\{1, \size{p}\}$. The function $\phi_L$ provides a bijection for all borderless pairs and almost half of the pairs whose border length set is $\{1, \size{p}\}$. These cases alone account for a sizable chunk of possible avoidant-equivalent word pairs, which is why the percentage of pairs for which we have natural bijections is so high. 
    
    \begin{table}[h!]
		\begin{tabular}{ |r||r||R{3cm}||r| }
			\hline
			Pattern length & $\phi_L$ bijection pairs & Composition bijection pairs& Equivalent pairs \\
			\hline
1&1&1&1\\
2&1&2&2\\
3&6&8&8\\
4&21&32&32\\
5&88&120&120\\
6&312&460&460\\
7&1212&1708&1716\\
8&4649&6764&6780\\
9&18264&26072&26168\\
10&71058&103460&103764\\
11&279946&403836&405404\\
12&1107836&1613132&1618556\\
			\hline
		\end{tabular} 
		\caption{Summary of bijections between patterns on $\{0,1\}$. }
		\label{bijections_table}
		\end{table}

 Table~\ref{bijections_table} contains data on the number of pairs of patterns on $\Sigma = \{0, 1\}$ for which we have a natural bijection. The ``Equivalent pairs'' column gives the total number of unordered pairs of patterns $p$ and $q$ for which $b(p)=b(q)$. The second column counts pairs for which $\phi_L$ establishes a bijection. Additionally, if we allow compositions with the reversal function and letter permutation function, we are able to obtain even more bijections; these pairs are counted in the third column.

	\section*{Acknowledgements}

    We thank the 2022 NYC Discrete Math REU for fostering a vibrant mathematical community rich in intellectual stimulation and gracious friendship alike. We extend our gratitude to the anonymous referee for the careful reading and valuable comments.
 
\bibliographystyle{plain}
\bibliography{refs}  

\begin{thebibliography}{1}

\bibitem{Claesson--Kitaev}
Anders Claesson and Sergey Kitaev.
\newblock Classification of bijections between 321- and 132-avoiding
  permutations.
\newblock {\em S\'{e}m. Lothar. Combin.}, 60:Art. B60d, 30, 2008/09.

\bibitem{Goulden--Jackson}
I.~P. Goulden and D.~M. Jackson.
\newblock An inversion theorem for cluster decompositions of sequences with
  distinguished subsequences.
\newblock {\em J. London Math. Soc. (2)}, 20(3):567--576, 1979.

\bibitem{Guibas--Odlyzko}
L.~J. Guibas and A.~M. Odlyzko.
\newblock String overlaps, pattern matching, and nontransitive games.
\newblock {\em J. Combin. Theory Ser. A}, 30(2):183--208, 1981.

\bibitem{holub}
\v{S}t\v{e}p\'{a}n Holub and Jeffrey Shallit.
\newblock Periods and borders of random words.
\newblock In {\em 33rd {S}ymposium on {T}heoretical {A}spects of {C}omputer
  {S}cience}, volume~47 of {\em LIPIcs. Leibniz Int. Proc. Inform.}, pages Art.
  No. 44, 10. Schloss Dagstuhl. Leibniz-Zent. Inform., Wadern, 2016.

\bibitem{Kim--Putcha--Roush}
Ki~Hang Kim, Mohan~S. Putcha, and Fred~W. Roush.
\newblock Some combinatorial properties of free semigroups.
\newblock {\em J. London Math. Soc.}, s2-16(3):397--402, 1977.

\bibitem{noonan}
John Noonan and Doron Zeilberger.
\newblock The {G}oulden--{J}ackson cluster method: extensions, applications and
  implementations.
\newblock {\em J. Differ. Equations Appl.}, 5(4-5):355--377, 1999.

\bibitem{silberger}
D.~M. Silberger.
\newblock How many unbordered words?
\newblock {\em Comment. Math. Prace Mat.}, 22(1):143--145, 1980.

\bibitem{solovev2}
A.~D. Solov’ev.
\newblock A combinatorial identity and its application to the problem
  concerning the first occurrence of a rare event.
\newblock {\em Theory of Probability \& Its Applications}, 11(2):276--282,
  1966.

\bibitem{Zeilberger}
Doron Zeilberger.
\newblock Enumeration of words by their number of mistakes.
\newblock {\em Discrete Mathematics}, 34(1):89--91, 1981.

\end{thebibliography}
 	
\end{document}